\documentclass[12pt]{article}

\title{Uniformly bounded arrays and mutually algebraic structures}

\author{ Michael C.\
Laskowski\thanks{Partially supported
by NSF grant DMS-1308546.}\\University of Maryland\and Caroline A.\ Terry\\University of Chicago}

\usepackage{amssymb}
\usepackage{amsthm}
\usepackage{amsmath}

\usepackage{url}
\usepackage{enumerate}
\usepackage{comment}

\usepackage{color}
\def\showauthornotes{1}

\ifnum\showauthornotes=1
\newcommand{\Authornote}[2]{{\sf\small\color{magenta}{[#1: #2]}}}
\else
\newcommand{\Authornote}[2]{}
\fi

\def\conc{{\char'136}}
\def\abar{\overline{a}}
\def\bbar{\overline{b}}
\def\cbar{\overline{c}}
\def\dbar{\overline{d}}
\def\ebar{\overline{e}}
\def\hbar{\overline{h}}

\def\rbar{\overline{r}}
\def\ubar{\overline{u}}

\def\wbar{\overline{w}}
\def\xbar{\overline{x}}
\def\ybar{\overline{y}}
\def\zbar{\overline{z}}

\def\phi{\varphi}

\def\B{{\cal B}}

\def\F{{\rm QF}}
\def\Fcal{{\cal F}}
\def\FF{{\bf F}}

\def\P{{\cal P}}

\def\S{{\cal S}}

\def\tp{{\rm tp}}

\def\acl{{\rm acl}}

\def\Fa0{{\FF^a_{\aleph_0}}}

\def\bp{{\bf Proof.}\quad}
\def\endproof{\medskip}
\def\<{\langle}
\def\>{\rangle}
\def\o2{{^{\omega} 2}}

\def\n2{{^{n} 2}}

\def\arr{{\rm arr}}

\def\QMA{{\rm QMA}}
\def\S{{\rm Supp}}
\def\M{{\mathbb M}}
\def\Z{{\mathbb Z}}
\def\AI{{\rm AI}}
\def\pp{{\bf p}}
\def\qq{{\bf q}}
\def\rr{{\bf r}}

\newtheorem{Theorem}{Theorem}[section]
\newtheorem{Proposition}[Theorem]{Proposition}
\newtheorem{Definition}[Theorem]{Definition}

\newtheorem{Remark}[Theorem]{Remark}

\newtheorem{Lemma}[Theorem]{Lemma}
\newtheorem{Corollary}[Theorem]{Corollary}
\newtheorem{Conclusion}[Theorem]{Conclusion}

\def\lg{{\rm lg}}

\newcommand\myrestriction{\mathord\restriction}
\def\mr#1{\myrestriction_{#1}}

\begin{document}

\date{\today}

\maketitle

\begin{abstract}  We define an easily verifiable notion of an atomic formula  having {\em uniformly bounded arrays} in a structure $M$.
We prove that if $T$ is a complete $L$-theory, then $T$ is mutually algebraic if and only if  there is some model $M$ of $T$ for which every atomic formula has uniformly bounded
arrays.  Moreover, an incomplete theory $T$ is mutually algebraic if and only if every atomic formula has uniformly bounded arrays in every model $M$ of $T$.
\end{abstract}

\section{Introduction}  
The notion of a mutually algebraic formula was introduced in \cite{DLR}, and the notions of  mutually algebraic structures and theories were introduced in \cite{MCL2}.  There, many properties were shown to be equivalent to 
mutual algebraicity, e.g.,  a structure $M$ is mutually algebraic if and only if every expansion $(M,A)$ by a unary predicate has the non-finite cover property (nfcp)
and a complete theory $T$ is mutually algebraic if and only if it is weakly minimal and trivial.  Whereas these characterizations indicate the strength of the hypothesis,
they do not lead to an easy verification that a specific structure is mutually algebraic.  The purpose of this paper is two-fold.  Primarily, we obtain equivalents of mutual algebraicity
that  are easily verifiable.  Most notably, we introduce the notion of a structure or a theory having {\em uniformly bounded arrays} and we
prove that for structures $M$ in a finite, relational language, $M$ is mutually algebraic if and only if $M$ has uniformly bounded arrays.  This result plays a key role in
\cite{MCLCAT}, where the authors describe the growth rates of hereditary classes of structures in a finite, relational language.  The other purpose is to develop a `Ryll-Nardzewski characterization' of mutual algebraicity, which is accomplished in Theorem~\ref{char}.  A theory $T$ is mutually algebraic if, for every restriction to a finite sublanguage, for every model $M$ and integer $n$, there are only finitely many quantifier-free $n$-types over $M$ that support an infinite array.

\section{Preliminaries}

Let $M$ be any $L$-structure and let $\phi(\zbar)$ be any $L(M)$-formula.  We say that $\phi(\zbar)$ is {\em mutually algebraic} if there is an integer $k$ such that
for any proper partition $\zbar=\xbar\conc\ybar$ (i.e., each of $\xbar,\ybar$ are nonempty) $M\models\forall\xbar\exists^{\le k}\ybar\phi(\xbar,\ybar)$.
Then, following \cite{MCL2}, a structure $M$ is {\em mutually algebraic} if every $L(M)$-formula is equivalent to a boolean combination of mutually algebraic $L(M)$-formulas,
and a theory $T$ is mutually algebraic if every model of $T$ is a mutually algebraic structure.  The following Theorem, which has the advantage of looking only at
atomic formulas, follows easily from two known results.  

\begin{Theorem}  \label{old}  Let $M$ be any $L$-structure.  Then $M$ is mutually algebraic if and only if every atomic formula $R(\zbar)$ is equivalent to a boolean combination of
quantifier-free mutually algebraic $L(M)$-formulas.
\end{Theorem}

\bp  First, assume $M$ is mutually algebraic.  The fact that every atomic $R(\zbar)$ is equivalent to a boolean combination of quantifier-free mutually algebraic $L(M)$-formulas
is the content of Proposition~4.1 of \cite{MCL1}.  For the converse, let $MA^*(M)$ denote the set of $L(M)$-formulas that are boolean combinations of mutually algebraic formulas.
This set is clearly closed under boolean combinations, and is closed under existential quantification by Propositon~2.7 of \cite{MCL2}.  Thus, if we assume that every atomic
formula is in $MA^*(M)$, it follows at once that every $L(M)$-formula is in $MA^*(M)$, hence $M$ is mutually algebraic.
\endproof

We will obtain a slight strengthening of Theorem~\ref{old} with Corollary~\ref{twopt}(2).
Whereas Theorem~\ref{old} placed no assumptions on the language, the main body of results in this paper assume that the underlying language is finite relational.
In Section~7 we obtain equivalents to mutual algebraicity for structures in arbitrary languages.

\medskip\noindent
{\bf Henceforth, for all results prior to Section~7, assume $L$ has finitely many relation symbols, finitely many constant symbols, and no function symbols.}

\medskip

For $L$ as above, fix an integer $n$ so that every atomic formula is at most $n$-ary.  As every quantifier-free type $p(\wbar)$ is determined by its family
of restrictions $\{p\mr{\wbar_i}$:  $\wbar_i$ a subsequence of $\wbar$ of length at most $n\}$, choose a specific $n$-tuple $\zbar=(z_1,\dots,z_n)$ of distinct
variable symbols.  Throughout, we concentrate on understanding spaces of quantifier-free types $p(\xbar)$, where $\xbar$ is a non-empty subsequence of $\zbar$.\footnote{The presentation of free variables in a type is  delicate, owing to the fact that `mutual algebraicity' is not preserved under adjunction of dummy variables.}

Fix an $L$-structure $M$.
For a non-empty subsequence $\xbar\subseteq\zbar$ and a subset $B\subseteq M$, let $\F_{\xbar}(B)$ denote the set of all quantifier-free
$L$-formulas $\phi(\xbar,\bbar)$ whose free variables are among $\xbar$ and $\bbar$ is from $B$.   Whereas we require $\xbar$ to be a subsequence of $\zbar$,
there are no limitations on the length of the parameter sequence $\bbar$.
By looking at the subsets of $M^{\lg(\xbar)}$ they
define, we can construe $\F_{\xbar}(B)$ as a boolean algebra.  Let $S_{\xbar}(B)$ denote its associated Stone space, i.e., the set of {\bf quantifier-free 
$\xbar$-types over $B$ that decide each $\phi\in\F_{\xbar}(B)$}.  
As usual, each of the Stone spaces $S_{\xbar}(B)$ are compact, Hausdorff, and totally disconnected when topologized by positing that the sets
$\{U_{\phi(\xbar,\bbar)}:\phi(\xbar,\bbar)\in\F_{\xbar}(B)\}$, where
$U_{\phi(\xbar,\bbar)}=\{p\in S_{\xbar}(B):\phi(\xbar,\bbar)\in p\}$, form a basis.
Moreover, because $L$ is finite relational, it follows  that  each $S_{\xbar}(B)$ is finite and every $p\in S_{\xbar}(B)$ is determined by a single $\phi(\xbar,\bbar)\in p$
whenever $B$ is finite.

\begin{Definition}  {\em  Fix a non-empty $\xbar\subseteq\zbar$, a subset $B\subseteq M$ and an integer $m$.  An $\xbar$-type $p\in S_{\xbar}(B)$ {\em supports an $m$-array} if
there is a {\bf pairwise disjoint} set $\{\dbar_i:i<m\}$ of (distinct) realizations of $p$ in $M$.  $p$ {\em supports an infinite array} if $M$ contains an infinite, pairwise disjoint set of realizations of $p$.
For each {\bf finite} $D\subseteq M$, let $N_{\xbar,m}(D)$ be the (finite) number of $p\in S_{\xbar}(D)$ that support an $m$-array.
}
\end{Definition}

The following definition is central to this paper, and forms the connection with \cite{MCLCAT}.  A local formulation, which relaxes the restriction on the language is given
in Section~7.

\begin{Definition} \label{restrict} {\em  A structure $M$ in a finite, relational language has {\em uniformly bounded arrays} if there is an integer $m>0$  such that
for every non-empty $\xbar\subseteq\zbar$, there is an integer $N$ such that $N_{\xbar,m}(D)\le N$
for all finite $D\subseteq M$.  When such an $N$ exists, we let $N_{\xbar,m}^{\arr}$ denote the smallest possible such $N$.
}
\end{Definition}  

It is easily seen that the properties described above are elementary.  In particular, if $m$ and the (finite) sequence $\<N_{\xbar,m}^\arr:\xbar\subseteq\zbar\>$ witness that
$M$ has uniformly bounded arrays, then the same $m$ and sequence $\<N_{\xbar,m}^\arr:\xbar\subseteq\zbar\>$ witness that
 any $M'$ elementarily equivalent to $M$ also has uniformly bounded arrays.  
 Because of this, we say that a complete theory $T$ in a finite, relational
 language {\em has uniformly bounded arrays} if some (equivalently all) models $M$ of $T$ have uniformly bounded arrays.

\section{Supportive and array isolating types}  
Throughout this and the next few  sections,  fix  a  complete theory $T$ in a finite, relational language $L$.  
Also fix an $\aleph_1$-saturated model $\M$ of $T$, which is a `monster model'
in the sense that all sets of parameters are chosen from $\M$.  The reader is reminded that all formulas and types mentioned are quantifier-free.

\begin{Definition}  {\em  For $\xbar$ a subsequence of $\zbar$ and $B$ countable, $\S_{\xbar}(B)$ is the set of all $p\in S_{\xbar}(B)$ that support an infinite array.
Let $\S(B)$ be the disjoint union of the spaces $\S_{\xbar}(B)$ for all subsequences $\xbar\subseteq\zbar$.
}
\end{Definition}

\begin{Lemma}  \label{compact}
Suppose $\xbar$ is a subsequence of $\zbar$, $B$ is countable, and $p\in S_{\xbar}(B)$.  
\begin{enumerate}
\item  The type $p\in\S_{\xbar}(B)$ if and only if for every $m\in\omega$ and every $\theta(\xbar,\bbar)\in p$, there is an $m$-array of solutions to $\theta(\xbar,\bbar)$.
\item  If $p$ has infinitely many solutions, then there is a (possibly empty) proper subsequence $\xbar_u\subseteq\xbar$ and a realization $\abar$ of $p$
such that $\tp((\abar\setminus\abar_u)/B\abar_u)$ supports an infinite array, where $\abar_u$ is the subsequence of $\abar$ corresponding to $\xbar_u$.
\end{enumerate}
\end{Lemma}

\bp  (1) is easily seen by compactness.  For (2), choose an infinite set $\{\abar_i:i\in\omega\}$ of distinct realizations of $p$.  If some infinite subset forms an array,
then take $\xbar_u=\emptyset$ and $p$ itself supports an infinite array.  If this is not the case, then
by the $\Delta$-system lemma,
there is an infinite $I\subseteq\omega$, a non-empty subsequence $\xbar_u$ of $\xbar$ and a fixed root $\rbar$ such that $(\abar_i)_u=\rbar$ and
$\abar_i\cap\abar_j=\rbar$ for distinct $i,j\in I$.  Then $\tp((\abar_i\setminus\rbar)/B\rbar)$ supports an infinite array whenever $i\in I$.
\endproof

It follows from Lemma~\ref{compact}(1)  that $\S_{\xbar}(B)$ is a closed, hence compact subspace of $S_{\xbar}(B)$.  
If we endow $\S(B)$ with the disjoint union topology (i.e., $U\subseteq\S(B)$ is open if and only if $(U\cap\S_{\xbar}(B))$ is open in $\S_{\xbar}(B)$ for every 
$\xbar\subseteq\zbar$) then $\S(B)$ is compact as well.

When our base set is a countable model, $\S_{\xbar}(M)$ is easily identified.

\begin{Lemma}  \label{genmodel}
If $M\preceq\M$ is countable, then $p\in S_{\xbar}(M)$ supports an infinite array if and only if $(x_i\neq a)\in p$ for all $x_i\in\xbar$ and all $a\in M$.
In particular, if $\cbar\cap M=\emptyset$, then $\tp(\cbar/M)\in \S_{\xbar}(M)$.
\end{Lemma}

\bp  Clearly, if $(x_i=a)\in p$ for any $x_i$ and any $a\in M$, then $p$ does not support a 2-array.
For the converse, choose any realization $\cbar$ of $p$ with $\cbar\cap M=\emptyset$.  To show that $p$ supports an infinite array, we employ Lemma~\ref{compact}(2).
Choose any $\theta(\xbar,\hbar)\in p$ (so $\hbar$ is from $M$).   We will construct an infinite array of solutions to $\theta(\xbar,\hbar)$ inside $M$.  
The construction is easy once we note that for any 
 finite subset $F\subseteq M$, $\cbar$ is a witness in $\M$ to
 $$\exists \xbar(\theta(\xbar,\hbar)\wedge\bigwedge_{x_i\in\xbar}\bigwedge_{a\in F} x_i\neq a)$$
 As $M\preceq\M$, we have some $\dbar\in M^{\lg(\xbar)}$ realizing $\theta(\xbar,\hbar)$ disjoint from $F$.
 \endproof

Next, we explore extensions of types $p\in\S(B)$.
By compactness, it is easily seen that whenever $B\subseteq B'$ are countable, then every $p\in \S(B)$ has an extension to some $q\in\S(B')$.
Abusing notation somewhat, let $\S(\M)$ denote the set of global types with the property that every restriction to a countable set supports an infinite array.
An easy compactness argument shows that every $p\in\S(B)$ has  a `global extension' to some $\pp\in \S(\M)$.
In general, a type $p\in\S(B)$ has many such global extensions, but we focus on when this is unique.

\begin{Definition}  {\em  A quantifier-free formula $\phi(\xbar,\ebar)$ is {\em array isolating} if there is exactly one global type $\pp\in\S_{\xbar}(\M)$
with $\phi(\xbar,\ebar)\in\pp$.   Call a global type $\pp\in\S_{\xbar}(\M)$ {\em array isolated} if it contains some array isolating formula.
Let $\AI_{\xbar}(\M)$ denote the set of array isolated global $\xbar$-types and let $\AI(\M)$ be the disjoint union of $\AI_{\xbar}(\M)$ over all
subsequences $\xbar\subseteq\zbar$.  For $\pp\in\AI(\M)$, $\pp|B$ denotes the restriction of $\pp$ to a type in $\S(B)$.
}
\end{Definition}

The following Lemma is immediate.  We will get a stronger conclusion in Section~\ref{52}  under the additional assumption that $\S_{\xbar}(\M)$ is finite.

\begin{Lemma}  \label{note}
Suppose $\pp\in\S_{\xbar}(\M)$  is array isolated as witnessed by the  array isolating $\phi(\xbar,\ebar)$.  Then any $q\in S_{\xbar}(B)$ containing $\phi(\xbar,\ebar)$
is either equal to $\pp|B$ or does not admit an infinite array.
\end{Lemma}

\bp  If $q\neq \pp|B$ supported an infinite array, then any global extension $\qq\supseteq q$ would be distinct from $\pp$.  This contradicts $\phi(\xbar,\ebar)$ being array isolating.
\endproof

\begin{Definition}  {\em  Say that $\pp\in\AI_{\xbar}(\M)$ is {\em based on $B$} if $\pp\cap\F_{\xbar}(B)$ contains an array isolating formula $\phi(\xbar,\ebar)$,
the interpretation $c^{\M}\in B$ for every constant symbol,
and,
moreover there is an infinite array $\{\abar_i:i\in\omega\}\subseteq B^{\lg(\xbar)}$ of pairwise disjoint realizations of $\phi(\xbar,\ebar)$.
}
\end{Definition}

Clearly, if $\pp$ is based on $B$, then it is also based on any $B'\supseteq B$.  If $B$ is a model, then the second and third clauses are redundant, that is:

\begin{Lemma} \label{inM}  If $M\preceq\M$, $\pp\in\AI_{\xbar}(\M)$  and $\pp\cap\F_{\xbar}(M)$ contains an array isolating formula $\phi(\xbar,\ebar)$,
then $\pp$ is based on $M$.
\end{Lemma}

\bp  As $M\preceq\M$, every $c^{\M}\in M$.  Now, fix an array  
isolating formula $\phi(\xbar,\ebar)\in\pp\cap\F_{\xbar}(M)$ and we recursively construct an infinite array of realizations of $\phi(\xbar,\ebar)$ inside $M$
as follows.  First, let $B_0=\ebar$  and let $p_0=\pp|B_0$.  As the language $L$ and $B_0$ are finite, $p_0$ is isolated by a formula over $B_0$.  As $M\preceq \M$,
choose a realization $\abar_0$ of $p_0$ inside $M$.  Then put $B_1=B_0\cup\{\abar_0\}$, let $p_1=\pp|B_1$, and continue for $\omega$ steps.
\endproof

There is a tight analogy between array isolated types $\pp$ based on $B$ and strong types  over $B$ in a stable theory, but in general they are not equivalent.
Indeed, as we are restricting to quantifier free types, a typical  restriction $\pp|B$ is not even a complete type with respect to formulas with quantifiers.
We  show that every $\pp\in\AI(\M)$ is  $B$-definable for any $B$ on which it is based.

\begin{Lemma}  \label{getm}  Suppose $\phi(\xbar,\ebar)$ is an array isolating formula and $\theta(\xbar,\ybar)\in\F_{\xbar\ybar}(\emptyset)$.
There is an integer $m=m(\phi(\xbar,\ebar),\theta)$ such that for all $\dbar\in \M^{\lg(\ybar)}$, exactly one of 
$\phi(\xbar,\ebar)\wedge\theta(\xbar,\dbar)$ and $\phi(\xbar,\ebar)\wedge\neg\theta(\xbar,\dbar)$ admits an $m$-array.
\end{Lemma}

\bp  As $\phi(\xbar,\ebar)$ admits an infinite array, at least one of the two formulas will as well.  However, if such an $m$ did not exist,
then for each $m$ there would be a tuple $\dbar_m$ such that both 
$\phi(\xbar,\ebar)\wedge\theta(\xbar,\dbar_m)$ and $\phi(\xbar,\ebar)\wedge\neg\theta(\xbar,\dbar_m)$ admit an $m$-array.
Thus, by the saturation of $\M$, there would be a tuple $\dbar^*$ such that both 
$\phi(\xbar,\ebar)\wedge\theta(\xbar,\dbar^*)$ and $\phi(\xbar,\ebar)\wedge\neg\theta(\xbar,\dbar^*)$ admit infinite arrays, contradicting $\phi(\xbar,\ebar)$ being
array isolating.
\endproof

\begin{Definition}  \label{definition}  {\em  Fix  any $\pp\in\AI_{\xbar}(\M)$ and any set $B$ on which it is based.
Choose an array isolating formula $\phi(\xbar,\ebar)\in\pp\cap\F_{\xbar}(B)$ and an infinite array
$\{\abar_i:i\in\omega\}\subseteq B^{\lg(\xbar)}$ of realizations of $\phi(\xbar,\ebar)$.
For any  $\theta(\xbar,\ybar)\in\F_{\xbar\ybar}(\emptyset)$  let
$$d_\pp\xbar \theta(\xbar,\ybar):=\bigvee_{s\in{{2m}\choose m}}\bigwedge_{i\in s} \theta(\abar_i,\ybar)$$
where $m=m(\phi(\xbar,\ebar),\theta)$ is chosen by Lemma~\ref{getm}.}
\end{Definition} 

Visibly, $d_\pp\xbar\theta(\xbar,\ybar)\in\F_{\ybar}(B)$.  Its relationship to $\theta(\xbar,\ybar)$  and $\pp$ is explained by the following Lemma.

  \begin{Lemma} \label{basedef}   Suppose $\pp\in\AI(\M)$  is based on  a countable set
  $B$ and $\phi(\xbar,\ebar)$ and $\{\abar_i:i\in\omega\}$ are chosen as in Definition~\ref{definition}.
The following are equivalent for any $\theta(\xbar,\ybar)\in\F_{\xbar\ybar}(\emptyset)$ and any $\dbar\in\M^{\lg(\ybar)}$:
\begin{enumerate}
\item $\M\models d_\pp\xbar\theta(\xbar,\dbar)$;
\item  $\theta(\xbar,\dbar)\in \pp$;
\item  For all countable $B'$, $\pp|B'\cup\{\theta(\xbar,\dbar)\}$ supports an infinite array; and
\item  The partial type $\pp|B\cup\{\theta(\xbar,\dbar)\}$ supports an array of length $m=m(\phi(\xbar,\ebar),\theta)$.
\end{enumerate}
\end{Lemma}

  \bp  
  $(1)\Rightarrow(2)$:  As $\M\models d_\pp\xbar\theta(\xbar,\dbar)$, some $m$-element subset of $\{\abar_i:i<2m\}$ is an $m$-array of realizations of 
  $\phi(\xbar,\ebar)\wedge\theta(\xbar,\dbar)$.  By choice of $m$, Lemma~\ref{getm} implies that $\phi(\xbar,\ebar)\wedge\theta(\xbar,\dbar)$ supports an
  infinite array, so 
  $\theta(\xbar,\dbar)\in \pp$.
  
  $(2)\Rightarrow(3)$:  Choose any countable $B'$.  If $\theta(\xbar,\dbar)\in \pp$, then as $\pp|B'\cup\{\theta(\xbar,\dbar)\}$ is a countable subset of $\pp$, 
it supports an infinite array.

  $(3)\Rightarrow(4)$:  Trivial.
  
  $(4)\Rightarrow(1)$: 
  Assume that  $\M\models \neg d_\pp\xbar\theta(\xbar,\dbar)$. Then some $m$-element subset of $\{\abar_i:i<2m\}$
  witnesses that  $\phi(\xbar,\ebar)\wedge\neg\theta(\xbar,\dbar)$ supports an $m$-array.  By Lemma \ref{getm}, $\phi(\xbar,\ebar)\wedge \theta(\xbar,\dbar)$ cannot support an $m$-array.  Consequently $p|B\cup \{\theta(\xbar,\dbar)\}$ cannot support an $m$-array (and thus cannot support an infinite array).
  \endproof

  \section{Free products of array isolated types}

Throughout this section, $T$ is a complete theory in a finite, relational language and $\M$ is an $\aleph_1$-saturated model, from which we take our parameters.  

In this section we describe how to construct a `free join' of array isolated types.
Suppose $\xbar,\ybar$ are disjoint, non-empty subsequences of $\zbar$,  $\pp(\xbar)\in\AI_{\xbar}(\M)$, and $\qq(\ybar)\in \AI_{\ybar}(\M)$.
We show that there is a well-defined $\rr(\xbar,\ybar)\in \S_{\xbar\ybar}(\M)$ constructed from this data.  We begin with lemmas that unpack our definitions.

\begin{Lemma}  \label{warm}
Suppose $\pp(\xbar),\qq(\ybar)$ are as above and $B$ is a countable set on which both $\pp$ and $\qq$ are based.
For any $\theta(\xbar,\ybar,\bbar)\in\F_{\xbar\ybar}(B)$ and any $\cbar$ realizing $\pp|B$,
$$\theta(\cbar,\ybar,\bbar)\in \qq|B\cbar \qquad \hbox{if and only if} \qquad  d_{\qq}\ybar\theta(\xbar,\ybar,\bbar)\in \pp|B$$
\end{Lemma}

\bp  First, assume $\theta(\cbar,\ybar,\bbar)\in \qq|B\cbar$.  Then $\qq|B\cup\{\theta(\cbar,\ybar,\bbar)\}\subseteq\qq$, hence it supports an infinite array.
By Lemma~\ref{basedef} applied to $\theta(\cbar,\ybar,\bbar)$ (i.e., taking $\dbar:=\cbar\bbar$), $$\M\models d_{\qq}\ybar\theta(\cbar,\ybar,\bbar)$$
Taking $\wbar$ to be a sequence of variables for $\bbar$, since $d_{\qq}\ybar\theta(\xbar,\ybar,\wbar)\in \F_{\xbar\wbar}(B)$, $\bbar$ is from $B$,
and $\cbar$ realizes $\pp|B$, we conclude that $d_{\qq}\ybar\theta(\xbar,\ybar,\bbar)\in \pp|B$.

The converse is dual, using $\neg \theta$ in place of $\theta$.
\endproof

\begin{Lemma}  \label{between}  Suppose $\pp(\xbar),\qq(\ybar)$ are as above and $B$ is a countable set on which both $\pp$ and $\qq$ are based.
For any $\theta(\xbar,\ybar,\bbar)\in\F_{\xbar\ybar}(B)$ and any $\cbar,\cbar'$ realizing $\pp|B$, 
$\theta(\cbar,\ybar,\bbar)\in\qq|B\cbar$ if and only if $\theta(\cbar',\ybar,\bbar)\in\qq|B\cbar'$.
\end{Lemma}

\bp  By Lemma~\ref{warm}, each statement is equivalent to $d_{\qq}\theta(\xbar,\ybar,\bbar)\in \pp|B$, which does not depend on our choice of $\cbar$.
\endproof

Extending this,

\begin{Lemma}  \label{same}
Suppose $\pp(\xbar),\qq(\ybar)$ are as above and $B$ is a countable set on which both $\pp$ and $\qq$ are based.
For any $\theta(\xbar,\ybar,\bbar)\in\F_{\xbar\ybar}(B)$ and any $\cbar,\cbar'$ realizing $\pp|B$, for any $\dbar$ realizing $\qq|B\cbar$ and $\dbar'$ realizing
$\qq|B\cbar'$, the following three notions are equivalent:
\begin{enumerate}
\item  $\M\models\theta(\cbar,\dbar,\bbar)$;
\item  $\M\models d_{\pp}\xbar [d_{\qq}\ybar\theta(\xbar,\ybar,\bbar)]$; and
\item  $\M\models\theta(\cbar',\dbar',\bbar)$.
\end{enumerate}
\end{Lemma}

\bp  Because of the duality in the statements, it suffices to prove $(1)\Leftrightarrow(2)$.
First, assume (1) holds.  As $\dbar$ realizes $\qq|B\cbar$, we infer $\theta(\cbar,\ybar,\bbar)\in \qq$.
So, by Lemma~\ref{basedef}, $\M\models d_{\qq}\ybar\theta(\cbar,\ybar,\bbar)$.
As $d_{\qq}\ybar\theta(\xbar,\ybar,\bbar)\in \F_{\xbar\ybar}(B)$ and $\cbar$ realizes $\pp|B$, we conclude that
$d_{\qq}\ybar\theta(\xbar,\ybar,\bbar)\in\pp$ and hence $\M\models d_{\pp}\xbar [d_{\qq}\ybar\theta(\xbar,\ybar,\bbar)]$.
Showing that $(\neg 1)$ implies $(\neg 2)$ is dual, using $\neg\theta$ in place of $\theta$.
\endproof

We now define the free product of array supporting global types.

\begin{Definition}  \label{free}  {\em
Suppose $\xbar,\ybar$ are disjoint subsequences of $\zbar$, $\pp\in\AI_{\xbar}(\M)$ and $\qq\in\AI_{\ybar}(\M)$.  
Then the {\em free product $\rr=\pp\times\qq$} is defined as
$$\rr(\xbar,\ybar):=\{\theta(\xbar,\ybar,\bbar)\in\F_{\xbar\ybar}(\M):\M\models d_{\pp}\xbar[d_{\qq}\ybar\theta(\xbar,\ybar,\bbar)]\}$$
}
\end{Definition}
Because of Lemma~\ref{same}, $\rr(\xbar,\ybar)$ is also equal to the set of all $\theta(\xbar,\ybar,\bbar)\in\F_{\xbar\ybar}(\M)$ such that  for some/every $B$ on which
both $\pp$ and $\qq$ are based and $\bbar$ is from $B$, 
for some/every $\cbar$ realizing $\pp|B$ and for some/every $\dbar$ realizing $\qq|B\cbar$ we have $\M\models\theta(\cbar,\dbar,\bbar)$.
It is easily seen from this characterization that $\rr(\xbar,\ybar)\in\S_{\xbar\ybar}(\M)$.

Next, we show that the free join is symmetric.  We begin with a Lemma.
\begin{Lemma} \label{1symm}
Suppose $\xbar,\ybar$ are disjoint subsequences of $\zbar$, $\pp(\xbar)\in\AI_{\xbar}(\M)$, $\qq(\ybar)\in\AI_{\ybar}(\M)$.
Then for every countable set $B$ on which both types are based, for every $\cbar$ realizing $\pp|B$, $\dbar$ realizing $\qq|B$, and for every 
$\theta(\xbar,\ybar,\bbar)\in\F_{\xbar\ybar}(B)$ such that $\M\models\theta(\cbar,\dbar,\bbar)$,
$$\theta(\xbar,\dbar,\bbar)\in \pp|B\dbar \qquad\hbox{if and only if}\qquad \theta(\cbar,\ybar,\bbar)\in \qq|B\cbar$$
\end{Lemma}

\bp  Assume by way of contradiction that $\theta(\xbar,\dbar,\bbar)\in \pp|B\dbar$, but  $\theta(\cbar,\ybar,\bbar)\not\in \qq|B\cbar$, with the other direction being dual.

Write $\theta$ as $\theta(\xbar,\ybar,\wbar)$.  As  both $\pp$ and $\qq$ are based on $B$, 
choose array isolating formulas $\phi(\xbar,\ebar)\in\F_{\xbar}(B)$ and $\psi(\ybar,\ebar')\in\F_{\ybar}(B)$ for $\pp$ and $\qq$, respectively.  
Let $m=\max\{m(\phi(\xbar,\ebar),\theta(\xbar;\ybar\wbar)), m(\psi(\ybar,\ebar'),\theta(\ybar;\xbar\wbar))\}$.  (Note the different partitions of $\theta$.)

As $B$ is countable, choose infinite arrays $\{\cbar_i:i\in\omega\}$ and $\{\dbar_j:j\in\omega\}$ for $\pp|B$ and $\qq|B$, respectively.  
By Lemma~\ref{warm} we have that $\theta(\xbar,\dbar_j,\bbar)\in \pp|B\dbar_j$ for each $j$ and that $\theta(\cbar_i,\ybar,\bbar)\not\in \qq|B\cbar_i$
for each $i$.  
By passing to infinite subsequences, we may additionally assume these sets are pairwise disjoint (i.e., $\cbar_i\cap\dbar_j=\emptyset$ for all $i,j$).
Choose a number $K>>m$.
Form a finite, bipartite graph with universe $C\cup D$, where $C=\{\cbar_i:i<K\}$ and $D=\{\dbar_j:j<K\}$ with an edge $E(\cbar_i,\dbar_j)$ if and only if
$\M\models\theta(\cbar_i,\dbar_j,\bbar)$.  We will obtain a contradiction by counting the number of edges in two different ways.

On one hand, because $\theta(\cbar_i,\ybar,\bbar)\not\in \qq|B\cbar_i$, 
$\qq|B\cup\{\theta(\cbar_i,\ybar,\bbar)\}$ does not support an infinite array for each $\cbar_i$.  By our choice of $m$,  it cannot support an array of length $m$ either.
Because any $m$-element subset of $D$ is an array of length $m$, we conclude that for every $\cbar_i$, there are fewer than $m$ many $\dbar_j$
such that $\M\models \theta(\cbar_i,\dbar_j)$.  Thus, the number of edges of the graph is bounded above by $Km$.
On the other hand, for any $\dbar_j$, as $\pp|B\cup\{\theta(\xbar,\dbar_j,\bbar)\}$ supports an infinite array, $\qq|B\cup\{\neg\theta(\xbar,\dbar_j,\bbar)\}$ cannot support an infinite
array, hence cannot support an array of size $m$.  Thus, the edge-valence of each $\dbar_j$ is at least $(K-m)$, implying that our graph has at least
$K(K-m)$ edges.  As $K$ is much larger than $m$, this is a contradiction.
\endproof

\begin{Corollary}  \label{2symm}
Suppose $\xbar,\ybar$ are disjoint subsequences of $\zbar$, $\pp\in\AI_{\xbar}(\M)$, and $\qq\in\AI_{\ybar}(\M)$.
Then $\pp\times\qq=\qq\times \pp$.  That is, for any set $B$ on which both $\pp,\qq$ are based and for any $\theta(\xbar,\ybar,\wbar)\in\F_{\xbar\ybar\wbar}(\emptyset)$,
the formulas $d_\pp\xbar [d_{\qq}\ybar \theta(\xbar,\ybar,\wbar)]$ and $d_\qq\ybar [d_{\pp}\xbar \theta(\xbar,\ybar,\wbar)]$  in $\F_{\wbar}(B)$ are equivalent.
\end{Corollary}

\bp  Choose any $\theta(\xbar,\ybar,\bbar)\in\pp\times\qq$ with $\bbar\in\M^{\lg(\wbar)}$. Choose any countable set $B$ containing $\bbar$ on which both $\pp$ and $\qq$ are based.
Choose $\cbar$ realizing $\pp|B$ and $\dbar$ realizing $\qq|B\cbar$.  By the equivalent definition of $\pp\times\qq$, $(\cbar,\dbar)$ realizes $\pp\times\qq$, hence
$\M\models\theta(\cbar,\dbar,\bbar)$.  Thus, by Lemma~\ref{1symm}, $\cbar$ also realizes $\pp|B\dbar$.  Hence $(\cbar,\dbar)$ also realizes $\qq\times\pp$, so
$\theta(\xbar,\ybar,\bbar)\in\qq\times\pp$ as well.
\endproof

\section{Finitely many mutually algebraic types supporting arrays}  \label{52}

We continue our assumption that $\M$ is an $\aleph_1$-saturated model of a complete theory $T$ in a finite, relational language.
We begin by recording two consequences of $\S_{\xbar}(\M)$ being finite.  Note that simply by adding repeated elements to tuples,
$\S_{\xbar}(\M)$ being finite implies $\S_{\xbar'}(\M)$ finite for all non-empty subsequences $\xbar'\subseteq\xbar$.

\begin{Lemma}  \label{suppfinite}  Fix $\xbar\subseteq\zbar$ and assume that $\S_{\xbar}(\M)$ is finite.  Then $\S_{\xbar}(\M)=\AI_{\xbar}(\M)$ and, moreover,
every $\pp\in\AI_{\xbar}(\M)$ is based on every $M\preceq\M$.
\end{Lemma}

\bp  As $\S_{\xbar}(\M)$ is always a closed subspace of $S_{\xbar}(\M)$, it is compact.  Thus, if it is finite, every $\pp\in\S_{\xbar}(\M)$ is isolated.  
For the moreover clause,  write $\S_{\xbar}(\M)=\{\pp_1,\dots,\pp_n\}$
and choose array isolating formulas $\phi_{i}(\xbar,\ebar_i)$ for each $\pp_i$.  By repeated use of Lemma~\ref{getm},
let $m^*$ be the maximum of all $m(\phi_i(\xbar,\ebar_i),R(\xbar,\ybar))$ among all $\pp_i\in\S_{\xbar}(\M)$ and all atomic $R\in L$.
Then $\M$ is a model of the sentence
\begin{quotation}
\noindent $\exists \wbar_1\dots\exists \wbar_n [\{\phi_i(\xbar,\wbar_i):1\le i\le n\}$ are pairwise inconsistent and,
 for all atomic $R\in L$, for all $\zbar$, and for all $1\le i\le n$, exactly one of
$\phi_i(\xbar,\wbar_i)\wedge R(\xbar,\zbar)$ and $\phi_i(\xbar,\wbar_i)\wedge\neg R(\xbar,\zbar)$ supports an $m^*$-array].
\end{quotation}
Thus, any $M\preceq\M$ also models this sentence.  Choose witnesses $\ebar_1^*,\dots \ebar_n^*$ from $M$.  Then, for each $i$,
there is a unique type in $S_{\xbar}(M)$ containing $\phi_i(\xbar,\ebar_i^*)$ and supporting an infinite array.
By a second use of $M\preceq\M$, each $\phi_i(\xbar,\ebar_i^*)$ array isolates a global type $\pp\in\S_{\xbar}(\M)$.  
As $|\S_{\xbar}(\M)|=n$, we conclude that every $\pp\in\S_{\xbar}(\M)$ is array isolated by some $L(M)$-formula.
 In light of Lemma~\ref{inM}, it follows that each $\pp$ is based on $M$.
\endproof

As a consequence of Lemma~\ref{suppfinite}, if $\S_{\xbar}(\M)$ is finite and $M\preceq\M$ is countable, then every $p\in S_{\xbar}(M)$ that supports an infinite array
contains an array isolating formula, hence has a unique supportive extension $\pp\in\S_{\xbar}(\M)$.  Thus, for any extension $q\in S_{\xbar}(M\cbar)$ of $p$,
either $q=\pp|M\cbar$ or else $q$ does not support an infinite array.  In fact, even more is true.

\begin{Definition} {\em  A type $q\in S_{\xbar}(M\cbar)$ {\em has a finite part} if there is some non-empty $\xbar'\subseteq\xbar$ and some $\theta(\xbar',\cbar,\hbar)\in q$
for which $\M\models \exists^{<\infty}\xbar' \theta(\xbar',\cbar,\hbar)$.
}
\end{Definition}

\begin{Lemma}  \label{finpart}
 Suppose $\S_{\xbar}(\M)$ is finite and $M\preceq\M$ is countable. 
For any $\pp\in\S_{\xbar}(\M)$,
if a type $q\in S_{\xbar}(M\cbar)$ extends $\pp|M$ but  $q\neq\pp|M\cbar$,
then $q$ has a finite part.
\end{Lemma}

\bp  We argue by induction on $\lg(\xbar)$, so assume that the statement holds for all proper $\xbar'\subsetneq\xbar$.  By way of contradiction, choose $\pp\in\S_{\xbar}(\M)$,
and $q\in S_{\xbar}(M\cbar)$ with  $\pp|M\subseteq q$,  $q\neq\pp|M\cbar$, but $q$ has no finite part. 
We will obtain a contradiction to Lemma~\ref{note} by showing that $q$ supports an infinite array.  
Toward that goal, choose any $\theta(\xbar,\cbar,\hbar)\in q$.  By Lemma~\ref{compact}(1), it suffices to find an infinite array of realizations to $\theta(\xbar,\cbar,\hbar)$.

As $q$ has no finite part, by taking $\xbar'=\xbar$, $q$ has infinitely many realizations.  Thus, by Lemma~\ref{compact}(2), there is 
 a proper subsequence $\xbar_u\subseteq\xbar$ and a realization $\abar$ of $q$ such that $\tp((\abar\setminus\abar_u)/M\cbar\abar_u)$ supports an infinite array.
 Trivially, if $\xbar_u$ is empty, then this type is $q$ itself, $q$ supports an infinite array.  So we assume $\xbar_u$ is non-empty.  Write
 $\xbar=\ybar\conc \xbar_u$ and for clarity, write $\rbar$ for $\abar_u$ and $\bbar$ for $(\abar\setminus\rbar)$, so $\abar=\bbar\conc\rbar$.
 
 Let $q^*(\ybar):=\tp(\bbar/M\cbar\rbar)$.  We know that  $\bbar'\conc\rbar$ realizes $q$ whenever $\bbar'$ realizes $q^*$.  
 Let $\pp^*(\ybar)\in\S_{\ybar}(\M)$ be any global type extending $q^*$.  As $\S_{\ybar}(\M)$ is finite, by Lemma~\ref{suppfinite} $\pp^*$ is based on $M$.
 Choose an array isolating formula $\phi(\ybar,\ebar)\in\pp^*$ with $\ebar$ from $M$, along with an infinite array $\{\bbar_i:i\in\omega\}\subseteq M^{\lg(\ybar)}$
 of realizations of $\phi(\ybar,\ebar)$.  As $\phi(\ybar,\ebar)\wedge\theta(\ybar,\rbar,\cbar,\hbar)\in q^*$,  Lemma~\ref{getm} implies that
 all but finitely many $\bbar_i$ realize $\theta(\ybar,\rbar,\cbar,\hbar)$, so by elimination and reindexing,
 assume they all do.  
 
 On the other hand, let
 $q_u(\xbar_u):=\tp(\rbar/M\cbar)$.  As $q$ does not have a finite part, neither does $q_u$.  As $\xbar_u$ is a proper subsequence of $\xbar$, our inductive hypothesis
 implies that $q_u$ must support an infinite array.  Let $\{\rbar_j:j\in\omega\}$ be such an array.  Note that as $\tp(\rbar_j/M\cbar)=\tp(\rbar/M\cbar)$ and 
 $\theta(\bbar_i,\xbar_u,\cbar,\hbar)\in q_u$, it follows that $\theta(\bbar_i,\rbar_j,\cbar,\hbar)$ holds for all $i,j\in\omega$.  From this, as both $\{\bbar_i:i\in\omega\}$
 and $\{\rbar_j:j\in\omega\}$ are infinite arrays, it is easy to meld subsequences of these to produce an infinite array $\{\bbar_k\rbar_k:k\in\omega\}$ of realizations of 
 $\theta(\ybar,\xbar_u,\cbar,\hbar)$.
 \endproof

Next, we add mutual algebraicity to the discussion of supportive and array isolating types.

\begin{Definition} {\em  For $\xbar\subseteq\zbar$ non-empty, a global, supportive type $\pp\in\S_{\xbar}(\M)$  is {\em quantifier-free mutually algebraic (QMA)} if 
$\pp$ contains a mutually algebraic formula $\phi(\xbar)\in\F_{\xbar}(\M)$. Let $\QMA_{\xbar}(\M)$ denote the set of QMA types in $\S_{\xbar}(\M)$.  
Let $\QMA(\M)$ be the (finite)
disjoint union of the sets $\QMA_{\xbar}(\M)$.  
}
\end{Definition}

The goal of this section will be to deduce consequences from $\QMA(\M)$ being finite.

\begin{Lemma} \label{53} If $\QMA_{\xbar}(\M)$ is finite, then every $\pp\in\QMA_{\xbar}(\M)$ is array isolated, i.e., $\pp\in\AI_{\xbar}(\M)$.
\end{Lemma}

\bp  Fix any $\pp\in\QMA_{\xbar}(\M)$.  Choose a mutually algebraic $\phi_0(\xbar,\ebar_0)\in\pp$.
For each $\qq\in\QMA_{\xbar}(\M)$ distinct from $\pp$, choose a formula $\phi_{\qq}(\xbar,\ebar_\qq)\in \pp\setminus\qq$.
Then the formula $\phi_0(\xbar,\ebar_0)\wedge\bigwedge_{\qq\neq\pp} \phi_{\qq}(\xbar,\ebar_\qq)$ array isolates $\pp$.
\endproof

\begin{Definition} \label{partition}  {\em  Fix any non-empty $\xbar\subseteq\zbar$.  A {\em partition  $\P=\{\xbar_1,\dots,\xbar_r\}$ of $\xbar$}
satisfies (1) each $\xbar_i$ non-empty and  (2)  Every $x\in\xbar$ is contained in exactly one $\xbar_i$.  
For $w\subseteq\{1,\dots,r\}$, let $\xbar_w$ be the subsequence of $\xbar$ with universe $\bigcup\{\xbar_i:i\in w\}$.

For $\cbar\in (\M)^{\lg(\xbar)}$, a partition $\P$ of $\xbar$ naturally induces a partition $\{\cbar_1,\dots,\cbar_r\}$ of $\cbar$.
For $w\subseteq \{1,\dots,r\}$, $\cbar_w$ is the subsequence of $\cbar$ corresponding to $\xbar_w$.
}
\end{Definition}

\begin{Definition}  {\em  Fix any $\xbar\subseteq\zbar$, any countable $M\preceq \M$,
and $\cbar\in (\M\setminus M)^{\lg(\xbar)}$.
A {\em maximal mutually algebraic decomposition of $\cbar$ over $M$} is a partition $\P=\{\xbar_1,\dots,\xbar_r\}$ of $\xbar$ for which the induced  partition
$\{\cbar_1,\dots,\cbar_r\}$ of $\cbar$ satisfies the following for each $i\in\{1,\dots,r\}$:
\begin{itemize}
\item  $\cbar_i$ realizes a mutually algebraic formula $\phi(\xbar_i)\in \F_{\xbar_i}(M)$; but 
\item For any proper extension $\xbar_i\subsetneq\ubar\subseteq\xbar$, the subsequence $\dbar$ of $\cbar$ induced by
$\ubar$ does not realize any mutually algebraic formula $\psi(\ubar)\in \F_{\ubar}(M)$.
\end{itemize}
}
\end{Definition}

\begin{Lemma}  \label{existence}  For any $\xbar\subseteq\zbar$ and every countable $M\preceq\M$,
 every $\cbar\in (\M\setminus M)^{\lg(\xbar)}$ admits a unique maximal mutually algebraic decomposition over $M$.
\end{Lemma}

\bp  First, by Lemma~\ref{genmodel}, both $\tp(\cbar/M)$ and $\tp(\cbar'/M)$ for any subsequence $\cbar'\subseteq\cbar$ support infinite arrays. 
Next, 
as every formula $\phi(x)$ in one free variable is mutually algebraic,  every singleton $c\in\cbar$ realizes a mutually algebraic formula.
 For each $x\in\xbar$, choose a subsequence $\xbar_i$ of $\xbar$ containing $x$ such that
$\cbar_i$ realizes a mutually algebraic formula in $\F_{\xbar_i}(M)$ and is {\em maximal} i.e., there is no proper extension $\xbar'\supsetneq\xbar_i$ for which
$\cbar'$ realizes a mutually algebraic formula in $\F_{\xbar'}(M)$.  Clearly, $\{\xbar_1,\dots,\xbar_r\}$ covers $\xbar$.
The fact that it is a partition follows from the fact that if $\xbar_i,\xbar_j$ are not disjoint and $\phi(\xbar_i)$, $\psi(\xbar_j)$ are each mutually algebraic,
then their conjunction $(\phi\wedge\psi)(\xbar_i\xbar_j)$ is mutually algebraic as well (see e.g. Lemma~2.4(6) of \cite{MCL1}).
\endproof

It is easily checked that if  $\{\cbar_1,\dots,\cbar_r\}$
is a maximal mutually algebraic decomposition of $\cbar$ over $M$,
then for any  $w\subseteq \{1,\dots,r\}$, the subset $\{\cbar_i:i\in w\}$ is a maximal, mutually algebraic decomposition of $\cbar_w$ over $M$.

We are now able to state and prove the following.

\begin{Proposition}  \label{bigind}
Suppose that $\QMA(\M)$ is finite.  Then, for every subsequence $\xbar\subseteq\zbar$, every
$\pp\in \S_{\xbar}(\M)$ is equal to a free product $\qq=\pp_1\times\dots\times\pp_r$ of types from $\QMA(\M)$.
In particular, each $\S_{\xbar}(\M)$ is finite.
\end{Proposition}

\bp  As the whole of $\QMA(\M)$ is finite, choose a finite $D$ such that for every $\xbar\subseteq\zbar$ and every $\qq\in\QMA_{\xbar}(\M)$,
there is a mutually algebraic formula $\gamma(\xbar)\in\qq\cap \F_{\xbar}(D)$.  Choose a countable $M\preceq\M$ with $D\subseteq M$.

We prove the Proposition by induction on $\xbar$, i.e., we assume the Proposition holds for all proper subsequences $\xbar'$ of $\xbar$ and prove the result for $\xbar$.
To base the induction, first note that if $x\in\zbar$ is a singleton, then  as every formula $\phi(x)$ is mutually algebraic, every $\qq\in\S_x(\M)$ is also in
$\QMA_x(\M)$, which we assumed was finite.

Now, suppose $\xbar$ is a subsequence of $\zbar$, $\lg(\xbar)\ge 2$, and the Proposition holds for every proper subsequence $\xbar'$ of $\xbar$.
In particular, as $\S_{\xbar'}(\M)$ is finite, Lemma~\ref{suppfinite} implies that each $\qq\in\S_{\xbar'}(\M)$  is based on $M$ and is also in $\AI(\M)$.

Choose any $\qq^*\in\S_{\xbar}(\M)$.  Towards showing that $\qq^*$ is a free product of types from $\QMA_{\xbar}(\M)$, choose any 
$\cbar\in(\M\setminus M)^{\lg{(\xbar)}}$ realizing $\qq^*|M$.
Suppose  the partition $\P=\{\xbar_1,\dots,\xbar_r\}$ of $\xbar$ yields the maximal, mutually algebraic decomposition
$\cbar_1\conc\dots\conc\cbar_r$ of $\cbar$ over $M$.  

There are now two cases.  First, if $r=1$, then $\tp(\cbar/M)$ contains a mutually algebraic formula, so 
$\qq^*\in\QMA_{\xbar}(\M)$ and we are finished.
So assume $r\ge 2$.  As notation, for each $1\le j\le r$, let $\wbar_j$ be the subsequence $\xbar_1\dots\xbar_{j-1}\xbar_{j+1}\dots\xbar_r$ of $\xbar$
and let $\dbar_j$ be the corresponding subsequence of $\cbar$.  As each $\dbar_j$ is a proper subsequence of $\cbar$, our inductive hypothesis implies that
$\tp(\dbar_j/M)$ contains an array isolating formula.
Let $\qq_j$ be the (unique) global extension of $\tp(\dbar_j/M)$ to $\AI_{\wbar_j}(\M)$.  
By our inductive hypothesis again, each $\qq_j=(\pp_1\times\dots \pp_{j-1}\times \pp_{j+1}\times \dots \pp_r)$.
As each $\pp_j$ is also array isolated, by iterating Lemma~\ref{1symm} finitely often it follows  that there is a unique supportive type $\rr^*(\xbar)$ which is equal to 
$\qq_j(\wbar_j)\times\pp_j(\xbar_j)$ for
every $1\le j\le r$.

In light of the characterization of free products following Definition~\ref{free},  in order to conclude that 
$\qq^*=\qq_r\times \pp_r=\rr^*$, it suffices to prove the following Claim.

\medskip
\noindent{\bf Claim.}  $\dbar_r$  realizes $\qq_r|M\cbar_r$.
\medskip

Assume this were not the case.
We obtain a contradiction by showing that the whole of $\tp(\cbar_1,\dots,\cbar_{r}/M)$ contains a mutually algebraic formula $\phi(\xbar)$.
As $\tp(\dbar_r/M)=\qq_r|M$, but $\tp(\dbar_r/M\cbar_r)\neq\qq_r|M\cbar_r$,  Lemma~\ref{finpart} allows us to
choose a maximal subsequence $\cbar_u$ of $\dbar_r$ and $\delta_r(\xbar_u,\cbar_{r},\bbar_r)\in\tp(\cbar_u/M\cbar_{r})$ ($\bbar_r$  from $M$)
with only finitely many solutions.  As $\tp(\cbar_i/M)$ is mutually algebraic for each $i$, it follows from the maximality of $u$ that there is a non-empty subset $w\subseteq\{1,\dots,{r-1}\}$
such that $\cbar_u=\cbar_w$.  To ease notation, say  $w=\{1,\dots, s\}$  for some $s\le r-1$.  
We argue that $s=r-1$.  If this were not the case, then  $\{\cbar_1,\dots,\cbar_s,\cbar_{r}\}$ would be
a maximal mutually algebraic decomposition over $M$  of the subsequence $\cbar_1\conc\dots\conc\cbar_s\conc\cbar_{r}$
whose corresponding variables $\xbar_1\dots\xbar_s\xbar_r$
form a proper subsequence of $\xbar$.    Thus, by our inductive hypothesis, $\tp(\cbar_1\dots\cbar_s\cbar_{r}/M)$
would equal $(\pp_1\times\dots\times \pp_s\times \pp_{r})|M$, which is contradicted by the formula 
$\delta_r(\xbar_1,\dots,\xbar_s,\xbar_{r},\bbar_r)\in \tp(\cbar_1\dots\cbar_s\cbar_{r}/M)$.
Thus, we conclude that $s=r-1$.  Hence, $\delta_r(\wbar_r,\cbar_r,\bbar_r)$ has only finitely many solutions.

Next, choose any $j< r$.  
Now the presence of the formula $\delta_r$ implies that $\qq^*\neq \rr^*$, hence $\qq^*\neq \qq_j\times \pp_j$.  From this, it follows
that $\tp(\dbar_j/M\cbar_j)\neq \qq_j|M\cbar_j$.  So, arguing just as above but replacing $r$ by $j$ throughout,
we conclude there is a formula $\delta_j(\xbar,\bbar_j)\in\tp(\cbar/M)$ for which $\delta_j(\wbar_j;\cbar_j\bbar_j)$ has only finitely many solutions.

Thus, if we choose a mutually algebraic formula $\gamma_j\in\tp(\cbar_j/D)$ for each $j\le r$, we conclude
that the formula 
$$\phi(\xbar_1,\dots,\xbar_{r},\bbar_1\dots\bbar_r)=
\bigwedge_{j\le r}\big(\gamma_j(\xbar_j) \wedge\delta_j(\wbar_j,\xbar_j,\bbar_j)\big)$$
is mutually algebraic with free variables $\xbar$ and is in $\tp(\cbar/M)$.  This contradicts our assumption that $\tp(\cbar/M)$ was not mutually algebraic.
This completes the proof of the Claim as well as the Proposition.
\endproof

\begin{Conclusion} \label{quote} If $\QMA(\M)$ is finite, then so is $\S(\M)$.  Moreover, for any $\xbar\subseteq\zbar$, any countable $M\preceq\M$,
 and any $\cbar\in(\M\setminus M)^{\lg(\xbar)}$,
 $\tp(\cbar/M)$ is determined by the maximal mutually algebraic partition $\P=\{\xbar_1,\dots,\xbar_r\}$ and the corresponding set  $\{\pp_1,\dots,\pp_r\}$ of  
$\QMA_{\xbar_i}(\M)$ types.
\end{Conclusion}

 \section{Mutual algebraicity and unbounded arrays}
 
 The whole of this section is devoted to the statement and proof of Theorem~\ref{char}.  
 It can be construed as a kind of `Ryll-Nardzewski theorem' for Stone spaces of quantifier-free types.

\begin{Theorem}  \label{char}  Suppose $T$ is a complete theory in a finite, relational language,
all of whose atomic formulas have free variables among $\zbar$, and let $\M$ be an $\aleph_1$-saturated model of $T$.
The following are equivalent.

\begin{enumerate}
\item  $T$ has uniformly bounded arrays;
\item  For all subsequences $\xbar\subseteq\zbar$, $\S_{\xbar}(\M)$ is finite;
\item  For all subsequences $\xbar\subseteq\zbar$, every global supportive type $\pp\in\S_{\xbar}(\M)$
is array isolated;
\item  Whenever $M\preceq N$ are models of $T$, for all subsequences $\xbar\subseteq\zbar$,
there are only finitely many types in $S_{\xbar}(M)$ realized in $(N\setminus M)^{\lg(\xbar)}$;
\item  For all models $M$ and all subsequences $\xbar\subseteq\zbar$, only finitely many types in $S_{\xbar}(M)$ both contain a mutually algebraic formula
and   support an infinite array; and
\item  $T$ is mutually algebraic.
\end{enumerate}
\end{Theorem}

\bp  We begin by showing that $(2)\Leftrightarrow(3)\Leftrightarrow(4)$.  Fix any non-empty subsequence $\xbar\subseteq\zbar$.  The key observation for showing
$(2)\Leftrightarrow(3)$ is that
if $X$ is any compact, Hausdorff space, then $X$ is finite if and only if every element $a\in X$ is isolated.
Suppose that (2) holds.  To establish (3), note that $\S_{\xbar}(\M)$ as a subspace of $S_{\xbar}(\M)$, the Stone space of all quantifier free types is closed, and hence
compact.  As (2) implies it is finite as well, every $\pp\in\S_{\xbar}(\M)$ must be isolated in the subspace, hence array isolated.

Verifying that $(3)\Rightarrow(2)$ uses the converse of this.  Fix $\xbar\subseteq\zbar$.  Applying (3) to the model $\M$ yields that
every element of $\S_{\xbar}(\M)$ is isolated.  As  $\S_{\xbar}(\M)$ is compact and Hausdorff, it must be  finite.

$(2)\Rightarrow(4)$ is easy.  Assume (2). It suffices to prove (4) for all countable $M\preceq N$.   As $\M$ is $\aleph_1$-saturated, we may assume $N\preceq\M$.
But now, by Lemma~\ref{genmodel}, for any $\cbar\in (N\setminus M)^{\lg(\xbar)}$, $\tp(\cbar/M)$ supports an infinite array.  As any $p\in\S_{\xbar}(M)$
extends to some $\pp\in\S(\M)$, there are only finitely many such types.

Next, suppose (2) fails.  Choose $\xbar\subseteq\zbar$ and a countable, infinite $Y\subseteq\S_{\xbar}(\M)$.  For each pair $\pp\neq\qq$, choose a formula
$\phi_{\pp\qq}(\xbar,\ebar_{\pp\qq})\in\pp\setminus\qq$.  Choose a countable $M\preceq\M$ containing $\{\ebar_{\pp\qq}:\pp\neq\qq\in Y\}$.
Thus, $\{\pp|M:\pp\in Y\}$ is a countably infinite set of types, each of which support an infinite array.  As $\M$ is $\aleph_1$-saturated, each such $\pp|M$
is realized by some  $\cbar\in\M^{\lg(\xbar)}$.  That $\cbar\cap M=\emptyset$ follows from the fact that $\pp|M$ supports an infinite array.

Continuing on, we consider $(2)\Rightarrow(5)$.  It suffices to prove this for $M$ countable, and we may
assume $M\preceq\M$. As any $p\in S_{\xbar}(M)$ supporting an infinite array has an extension to $\S_{\xbar}(\M)$, (5) follows from (2).

$(5)\Rightarrow(2)$ is immediate from Conclusion~\ref{quote}.

$(2)\Rightarrow(1)$.  By (2), for each $\xbar\subseteq\zbar$, let $N(\xbar):= |\S_{\xbar}(\M)|$.    As $(2)\Rightarrow(3)$, every $p\in\S(\M)$
is array isolated.  For each $\xbar\subseteq\zbar$ and each $p\in\S_{\xbar}(\M)$, choose an array isolating formula $\phi_\pp(\xbar,\ebar_\pp)\in\pp$ and let $D\subseteq\M$ be finite and contain all $\ebar_\pp$ 
for all $\pp\in\S(\M)$.  It is easily seen that for every $\xbar\subseteq\zbar$ and every countable $D\subseteq B\subseteq\M$, $S_{\xbar}(B)$ has exactly $N(\xbar)$ types that support an infinite array.  Moreover, each such $q\in\S_{\xbar}(B)$ has a unique restriction $q|D\in\S_{\xbar}(D)$ and a unique extension $\qq\in\S_{\xbar}(\M)$.

Towards finding an appropriate $m$ as in Definition~\ref{restrict}, fix $\xbar\subseteq\zbar$ and partition each atomic $R(\zbar)\in L$ as $R(\xbar,\wbar)$.  
For each $\pp\in\S_{\xbar}(\M)$ with array isolating formula $\phi_\pp(\xbar,\ebar_\pp)\in\pp$, let $m(\pp,\xbar)$ be the maximum of the $2|L|$ numbers 
$m(\phi_\pp(\xbar,\ebar_\pp), \pm R(\xbar,\wbar))$
obtained by Lemma~\ref{getm}. That is, apply the Lemma $2|L|$ times, once for each $R\in L$, and once for each $\neg R$ for $R\in L$.

The point is that if $B$ is countable, $D\subseteq B\subseteq\M$, and $q\in S_{\xbar}(B)$ contains some $\phi(\xbar,\ebar_\pp)$ and supports an $m(\pp,\xbar)$-array,
then for every $R(\xbar,\bbar)$, $R(\xbar,\bbar)\in q$ if and only if $\phi_\pp(\xbar,\ebar_\pp)\wedge R(\xbar,\bbar)$ supports an $m(\pp,\xbar)$-array if and only if
$\phi_\pp(\xbar,\ebar_\pp)\wedge R(\xbar,\bbar)$ supports an infinite array if and only if $R(\xbar,\bbar)\in \pp$.
Thus, $q\subseteq \pp$.

On the other hand, let $\theta(\xbar):=\bigwedge\{\neg\phi_{\pp}(\xbar,\ebar_{\pp}):\pp\in \S_{\xbar}(\M)\}$.  Since there is no $q\in S_{\xbar}(D)$  with
$\theta\in q$ that supports an infinite array, compactness yields an integer $m^*(\xbar)$ such that no $q\in S_{\xbar}(D)$ with $\theta\in q$ supports an
$m^*(\xbar)$-array.  Clearly, for any $B\supseteq D$, no $q\in S_{\xbar}(B)$ with $\theta\in q$ could support an $m^*(\xbar)$-array either.  

Choose an integer $m$ that is greater than all $m^*(\xbar)$ and all $m(\pp,\xbar)$ for  $\xbar\subseteq\zbar$ and $\pp\in\S_{\xbar}(\M)$.
Combining the statements above, we see that for any countable $B\supseteq D$ and any $\xbar\subseteq\zbar$,
exactly $N(\xbar)$ types in $S_{\xbar}(B)$ support $m$-arrays.  Thus, $\M$ (and hence $T$ by elementarity) has uniformly bounded arrays.

$(1)\Rightarrow(2)$ is also easy.  Assume $T$ has uniformly bounded arrays.  Choose 
$m$ and $\<N^\arr_{\xbar,m}:\xbar\subseteq\zbar\>$ from the definition.  To establish (2), we claim that $|\S_{\xbar}(\M)|\le N^\arr_{\xbar,m}$ for each $\xbar\subseteq\zbar$.
To see this, fix $\xbar\subseteq\zbar$ and assume by way of contradiction there is a finite $Y\subseteq\S_{\xbar}(\M)$ with $|Y|>N^\arr_{\xbar,m}$.
For each pair $\pp\neq\qq$ from $Y$, choose some $\phi_{\pp\qq}(\xbar,\ebar_{\pp\qq})\in\pp\setminus \qq$ and let $D\subseteq\M$ be finite, containing all these
$\ebar_{\pp\qq}$.  Thus, the set $\{\pp|D:\pp\in Y\}$ are distinct elements of $S_{\xbar}(D)$.  As each restriction $\pp|D$ supports an infinite array, each supports an
$m$-array, contradicting our definition of $N^\arr_{\xbar,m}$.

Thus, conditions (1)--(5) are equivalent.

$(6)\Rightarrow(5)$: 
Suppose $T$ is mutually algebraic.  Then by Theorem~3.3  of \cite{MCL2},
 $T$ is weakly minimal, trivial, has nfcp, and moreover, any expansion of any model of $T$ by unary predicates also has an nfcp theory.  
 By way of contradiction, fix $M$, $\xbar\subseteq\zbar$, and an infinite set of distinct mutually algebraic types $\{p_i:i\in\omega\}\subseteq\S_{\xbar}(M)$.  
 As there are only finitely many permutations of $\xbar$, we may assume that for distinct $i,j\in\omega$, no permutation of a realization of $p_i$ realizes $p_j$.
 As the language is finite relational, choose the set $\Fcal$ of $L$-formulas as in Proposition~\ref{ap}.  By the pigeon-hole principle and relabelling,
 choose an $L$-formula $\phi(\xbar,\wbar)\in\Fcal$ such that for every $i\in\omega$, there is $\hbar_i\in M^{\lg(\wbar)}$ such that $\phi(\xbar,\hbar_i)\in p_i$
 and $\phi(\xbar,\hbar_i)$ is mutually algebraic.  Note that since $T$ has nfcp,  there is an $M$-definable formula $\mu(\wbar)$ such that for any $\hbar\in M^{\lg(\wbar)}$,
 $M\models\mu(\hbar)$ if and only if $\phi(\xbar,\ybar)$ is mutually algebraic.

 Let $\M\succeq M$ be $|M|^+$-saturated.  We will obtain a contradiction by constructing an expansion $\M^+=(\M,U,V,W)$ by unary predicates that has the finite cover property.
 From the saturation of $\M$, choose $k$-tuples $\{\abar^j_i:i\in\omega,j\le i\}$ (where $k=\lg(\xbar)$) such that $\abar_i^j$ realizes $p_i(\xbar)$ but
 $\acl(M\abar_i^j)\cap \acl(M\abar_{i'}^{j'})=M$ unless $i=i'$ and $j=j'$.  Let $A=\bigcup\bigcup\{\abar_i^j:i\in\omega,j\le i\}$ and let $B=\{b_i^j:i\in\omega,j\le i\}$
 be the subset of $A$ consisting of the `first coordinates' i.e., $b_i^j=(\abar_i^j)_0$ for all $i,j$.  

 Let $\M^+$ be the expansion of $\M$ by interpreting $U^{\M^+}=A$, $V^{\M^+}=B$, and $W^{\M^+}=M$.
 Let $$P(\xbar):=\bigwedge_{x\in\xbar} U(x)\wedge \exists\wbar\,\big[\bigwedge_{w\in\wbar} W(w)\wedge \mu(\wbar)\wedge \phi(\xbar,\wbar)\big]$$
 Note that $\M^+\models P(\abar_i^j)$ for all $i\in\omega$, $j\le i$.
 Also, if  $b_i^j\in B$  and $\M^+\models P(b_i^j,a_2,\dots a_k)$, then by mutual algebraicity, $\{a_2,\dots,a_k\}\subseteq \acl(M\abar_i^j)$, so $(b_i^j,a_2,\dots,a_k)$ is a permutation
 of $\abar_i^j$.
Let $\xbar'=(x_2,\dots,x_k)$, $\ybar'=(y_2,\dots,y_k)$,  and put
$$E(x,y):=\exists \xbar'\exists \ybar'[P(x,\xbar')\wedge P(y,\ybar')\wedge (\forall\wbar\in W)\bigwedge_{R\in L} R(x,\xbar',\wbar)\leftrightarrow R(y,\ybar',\wbar)]$$
$E$
is an $\M^+$-definable equivalence relation on $B=V^{\M^+}$, and $\M^+\models E(b_i^j,b_{i'}^{j'})$ if and only if
$i=i'$.   Thus, $E$ has arbitrarily large finite classes, which contradicts $Th(\M^+)$ having nfcp.

 \begin{Remark}  {\em  The implication $(6)\Rightarrow(5)$ above really relies on counting {\em quantifier-free} mutually algebraic types that support infinite arrays.  As an example,
 $Th(\Z,S)$ is mutually algebraic, but there are infinitely many mutually algebraic formulas $\phi_n(x,y)$, each of which support infinite arrays.
 Take $\phi_n(x,y):=\exists z_0\dots\exists z_n [(x=z_0)\wedge (y=z_n) \wedge\bigwedge_{i=0}^{n-1} S(z_n,z_{n+1})]$.
 }
 \end{Remark}

$(\hbox{2--5})\Rightarrow(6)$:   We assume all of (2--5) and prove that $T$ is mutually algebraic by way of Theorem~\ref{old}.  Choose an $\aleph_1$-saturated $\M\models T$.
As $\QMA(\M)$ is finite, choose a finite $D\subseteq\M$ so that every $\pp\in\QMA(\M)$ contains a mutually algebraic formula in $\F(D)$, and choose a countable $M\preceq\M$
with $D\subseteq M$.
 Fix $\xbar\subseteq\zbar$ and let
\begin{quotation}
\noindent $\B_{\xbar}:=\{$all  $L(M)$-formulas $\phi(\xbar)$ that are equivalent to boolean combinations of {\bf quantifier-free} mutually algebraic $\alpha(\xbar')$ for some $\xbar'\subseteq\xbar\}$.
\end{quotation}
Note that since $y=y$ is mutually algebraic and $\phi(\xbar):=\alpha(\xbar')\wedge\bigwedge_{y\in\xbar\setminus\xbar'} y=y$ is equivalent to $\alpha(\xbar')$, we can construe
every quantifier-free, mutually algebraic $\alpha(\xbar')$ as an element of $\B_{\xbar}$.

\medskip\noindent
{\bf Claim 1.}  If $\cbar,\cbar'\in \M^{\lg(\xbar)}$ satisfy 
$\phi(\cbar)\leftrightarrow\phi(\cbar')$ for every $\phi\in\B_{\xbar}$, then $\tp(\cbar/M)=\tp(\cbar'/M)$.

\bp  First, as $x=m$ is mutually algebraic for each $x\in\xbar$ and $m\in M$, it suffices to show this for all
$\cbar,\cbar'\in(\M\setminus M)^{\lg(\xbar)}$.  In this case, as $\cbar,\cbar'$ agree on all quantifier-free, mutually algebraic $L(M)$-formulas
$\alpha(\xbar_u)$ for $\xbar_u\subseteq\xbar$, it follows that $\cbar=\cbar_1\dots\cbar_r$ and $\cbar'=\cbar_1'\dots\cbar_r'$ have corresponding
maximal mutually algebraic decompositions over $M$, and moreover, by our choice of $D$, $\tp(\cbar_i/M)=\pp|M$ if and only if $\tp(\cbar_i'/M)=\pp|M$
for every $1\le i\le r$ and every $\pp\in\QMA(\M)$.  Thus, $\tp(\cbar/M)=\tp(\cbar'/M)$ by Conclusion~\ref{quote}.
\endproof

\medskip\noindent
{\bf Claim 2.}  Every $\theta(\xbar)\in\F_{\xbar}(M)$ is equivalent to some $\phi(\xbar)\in\B_{\xbar}$.

\bp  Fix $\theta(\xbar)\in\F_{\xbar}(M)$.  We first show that for every $\ebar\in\theta(\M)$, there is some $\delta_{\ebar}(\xbar)\in\B_{\xbar}$
such that $\M\models \delta_{\ebar}(\ebar)$ and $$\M\models\forall\xbar (\delta_{\ebar}(\xbar)\rightarrow\theta(\xbar)).$$
To see this, let $\Gamma(\xbar):=\{\delta(\xbar)\in\B_{\xbar}:\M\models\delta(\ebar)\}$.  By Claim 1, every $\ebar'$ realizing $\Gamma(\xbar)$ also realizes
$\theta(\xbar)$, so the existence of $\delta_{\ebar}(\xbar)$ follows by compactness and the $\aleph_1$-saturation of $\M$.
Let $\Delta(\xbar)=\{\delta_{\ebar}(\xbar):\ebar\in\theta(\M)\}$.  A second application of compactness and saturation implies that
some finite $\Delta_0\subseteq\Delta$ satisfies 
$$\M\models \forall\xbar(\theta(\xbar)\rightarrow\bigvee\Delta_0(\xbar))$$
so take $\phi(\xbar):=\bigvee\Delta_0(\xbar)$.
\endproof

That $T=Th(M)$ is mutually algebraic now follows immediately from Claim 2 and Theorem~\ref{old} applied to $M$.
\endproof

\section{Identifying mutually algebraic structures and theories in arbitrary languages}

In this section, we use Theorem~\ref{char} to deduce  local tests for mutual algebraicity without regard to the size of the language, nor the completeness of the theory.
We begin by noting, either by Lemma~2.10 of \cite{MCL2} or deducing it from Lemma~\ref{old}, that `mutual algebraicity' is an elementary property, i.e., if
$M$ is mutually algebraic, then any elementarily equivalent $M'$ is mutually algebraic as well.

\begin{Lemma}  \label{shrink}  Suppose $L_0\subseteq L$ are finite relational languages, $M$ is an $L$-structure, and $M_0$ is its reduct to an $L_0$-structure.
\begin{enumerate}
\item  If $D\subseteq M$ is finite, every $L_0$-type $p\in\S_{\xbar}^{L_0}(D)$ has an extension to an $L$-type $q\supseteq p$ with $q\in \S^{L}_{\xbar}(D)$.
\item  If $M$ has uniformly bounded arrays, then so does $M_0$.
\end{enumerate}
\end{Lemma}

\bp (1)  Choose an infinite array $\{\abar_i:i\in\omega\}\subseteq (M_0)^{\lg(\xbar)}$ of realizations of $p$.  As both $D$ and $L$ are finite, there are only finitely
many $L$-types in $S_{\xbar}(D)$, so there is an infinite
subsequence $\{\abar_i:i\in I\}$ such that $\tp_M(\abar_i/D)=\tp_M(\abar_j/D)$ for all $i,j\in I$.  Then $\tp_M(\abar_i)$ for any $i\in I$ is as required.

(2)  Assume that $M$ has uniformly bounded arrays.
Let $\M$ be an $\aleph_1$-saturated model of $Th(M)$ and let $\M_0$ be the reduct of $\M$ to $L_0$.  So $\M_0$ is an $\aleph_1$-saturated model of $Th(M_0)$.
Suppose every $L$-atomic formula has free variables among $\zbar$. 
By applying Theorem~\ref{char} to $Th(M)$, $\S^L_{\xbar}(\M)$ is finite for all subsequences $\xbar\subseteq\zbar$, and by a second application of Theorem~\ref{char}
to $Th(M_0)$, it suffices to establish the following Claim:

\medskip\noindent{\bf Claim.}  For each $\xbar\subseteq\zbar$, $|\S_{\xbar}^{L_0}(\M_0)|\le |\S_{\xbar}^{L}(\M)|$.

\medskip\bp   Fix $\xbar\subseteq\zbar$, and assume that $\S_{\xbar}^{L}(\M)=\{\qq_i:i<N\}$.
In order to establish the Claim, it suffices to prove that for every finite $D\subseteq\M_0$, every $L_0$-type $p\in\S^{L_0}_{\xbar}(D)$ is a subset of some $\qq_i$.
To see this, choose
any finite $D\subseteq \M_0$ and any $p\in \S^{L_0}_{\xbar}(D)$.  By (1), there is some $L$-type $q\supseteq p$ with $q\in\S^L_{\xbar}(D)$.  Since $q$ has an extension to
some $\qq_i\in\S^{L}_{\xbar}(\M)$, it follows that $p\subseteq\qq_i$ for some $i<N$.
\endproof

\begin{Definition}  \label{R}  {\em  For $L$ an arbitrary language, fix any atomic $L$-formula $R(\zbar)$.   Let $L_R:=\{R,=\}$, which is visibly finite relational.  
Given an $L$-structure $M$, let $M_R$ denote the reduct of $M$ to an $L_R$-structure.  We say {\em $R$ has uniformly bounded arrays in $M$} if
the $L_R$-structure $M_R$ has uniformly bounded arrays.
}
\end{Definition}

\begin{Theorem} \label{local}  The following are equivalent for an  $L$-structure $M$ in an arbitrary  language $L$:
\begin{enumerate}
\item  $M$ is mutually algebraic;
\item  Every  atomic $R(\zbar)$  has uniformly bounded arrays in $M$; 
\item  For every atomic $R(\zbar)$, the reduct $M_R$ is mutually algebraic.
\end{enumerate}
\end{Theorem}

\bp  The equivalence of (2) and (3) follows by applying Theorem~\ref{char} to each of the (finite relational) $L_R$-theories $Th(M_R)$.

For $(3)\Rightarrow(1)$, in order to prove that $M$ is mutually algebraic, by Theorem~\ref{old}, it suffices to prove that every atomic $R(\zbar)$ is equivalent to
a boolean combination of mutually algebraic formulas.  Fix an atomic $R(\zbar)$.  By (3) and Theorem~\ref{old} applied to $M_R$, $R(\zbar)$ is equivalent to a boolean
combination of quantifier-free mutually algebraic $L_R(M)$-formulas.  As a mutually algebraic formula in $M_R$ is also mutually algebraic in $M$, the result follows.

Finally, assume (1).  To obtain (2), fix an atomic $R(\zbar)$.  By Theorem~\ref{old}, choose a finite set 
$\{\phi_1(\xbar_1,\ebar_1),\dots,\phi_k(\xbar_k,\ebar_k)\}$ 
of mutually algebraic, quantifier-free $L$-formulas for which $R(\zbar)$ is equivalent in $M$ to some boolean combination (so each $\xbar_i$ is a subsequence of $\zbar$).
Expand $L$ to $L'$, adding new $\lg(\xbar_i)$-ary relation symbols $U_i$ and let $M'$ be the definitional expansion interpreting each $U_i$ as $\phi_i(M,\ebar_i)$.
Let  $L_0=\{U_1,\dots,U_k\}$, $L_0^R=L_0\cup\{R\}$, and let $M_0, M_0^R$ be the reducts of $M'$ to $L_0$ and $L_0^R$, respectively.  
Note that $M_0$ and $M_0^R$ have the same quantifier-free definable sets, and that the reduct of  $M_0^R$ to $L_R$ is $M_R$.

As each $L_0$-atomic formula is mutually algebraic, it follows from
Theorem~2.7 of \cite{MCL2} that $M_0$ is mutually algebraic.  
As $L_0$ is finite relational, by applying Theorem~\ref{char} to $Th(M_0)$, $M_0$ has uniformly bounded arrays.  
Since  $M_0$ and $M_0^R$ have the same quantifier-free definable sets, we conclude that $M_0^R$ also has uniformly bounded arrays.
As $L_0^R$ is finite relational, it follows from Lemma~\ref{shrink} that $M_R$  has uniformly bounded arrays as well.
\endproof

The following Corollary now follows easily.  Clause~2 is a slight strengthening of Theorem~\ref{old}.

\begin{Corollary}  \label{twopt}  Let $L$ be an arbitrary language.
\begin{enumerate}
\item  The reduct of a mutually algebraic $L$-structure is mutually algebraic.
\item  If $M$ is mutually algebraic, then every atomic $R(\zbar)$ is equivalent to a boolean combination of mutually algebraic, quantifier-free
$L_R$-formulas.
\end{enumerate}
\end{Corollary}

\bp  (1)  Let $M$ be any mutually algebraic $L$-structure, let $L_0\subseteq L$ be arbitrary, and let $M_0$ be the reduct of $M$ to $L_0$.
Fix any atomic $R(\zbar)\in L_0$.  
Applying  Theorem~\ref{local} to $M$ gives $M_R$ mutually algebraic.  As this holds for all atomic $R\in L_0$,  a second application of
Theorem~\ref{local} implies $M_0$ is mutually algebraic.

(2)  is also by Theorem~\ref{local}.
\endproof

Finally, we consider incomplete theories.  The following Corollary follows immediately from Theorem~\ref{local}, 
as by definition, an incomplete theory $T$ is mutually algebraic if and only if
every model $M\models T$ is mutually algebraic.

\begin{Corollary}  \label{inc}
 A possibly incomplete theory $T$  in an arbitrary  language  
 is mutually algebraic if and only if for every $M\models T$, every atomic $R(\zbar)$ has uniformly bounded arrays in $M$.
\end{Corollary}

\appendix

\section{A basis of mutually algebraic formulas}

\begin{Lemma}  \label{conj}  
Let $T$ be an arbitrary $L$-theory, $\lg(\xbar)=k$, and suppose that 
$$\theta(\xbar):=\bigwedge_{i\in S} \alpha_i(\xbar_i)\wedge\bigwedge_{j\in U} \neg \beta_j(\xbar_j)$$ is mutually algebraic and supports an infinite array, with each $\alpha_i(\xbar_i),\beta(\xbar_j)$
mutually algebraic and each $\xbar_i,\xbar_j$ a subsequence of $\xbar$.
Then there is a subset $S_0\subseteq S$ of size at most $k$ such that $\bigcup_{i\in S_0}\xbar_i =\xbar$ and $\theta^*(\xbar)=\bigwedge_{i\in S_0} \alpha_i(\xbar_i)$ is mutually algebraic.
\end{Lemma}

\bp Form a maximal sequence $\<i_0,\dots,i_{m-1}\>$ from $S$ such that 
for each $j<m$, $\xbar_{i_j}$ is properly partitioned by $\bigcup_{t<j}\xbar_{i_t}$, i.e.,  $\xbar_{i_j}\cap\bigcup_{t<j}\xbar_{i_t}\neq\emptyset$ and 
$\xbar_{i_j}\setminus\bigcup_{t<j}\xbar_{i_t}\neq\emptyset$. As $\lg(\xbar)=k$, $m\le k$.  Take $S_0=\{i_j:j<m\}$, let $\xbar_m=\bigcup_{i\in S_0} \xbar_i$.
By iterating Lemma~2.4(6) of \cite{MCL1},
 $\theta^*(\xbar_m):=\bigwedge_{i\in S_0} \alpha_i(\xbar_i)$ is mutually algebraic, so to complete the proof it suffices to prove that $\xbar_m=\xbar$.
 
 Suppose this were not the case, i.e., write $\xbar=\xbar_m\conc\ybar$ with $\ybar\neq\emptyset$.  
Choose an infinite array $\{\abar_n:n\in\omega\}$ of realizations of $\theta(\xbar)$.    Let $\bbar=\abar_0\mr{\ybar}$ and, for every $n\in\omega$, let
$\cbar_n=(\abar_n\mr{\xbar_m})\conc\bbar$.  It suffices to prove the following Claim, as it contradicts $\theta(\xbar)$ being mutually algebraic.

\medskip\noindent{\bf Claim.}  $\theta(\cbar_n)$ holds for cofinitely many $n$.

\medskip  Towards the Claim, we first show that $\alpha_i(\cbar_n\mr{\xbar_i})$ holds for every $i\in S$ and $n\in\omega$.
To see this, fix $i\in S$.  By the maximality of the sequence defining $S_0$, either $\xbar_i\subseteq\xbar_m$ or $\xbar_i\cap\xbar_m=\emptyset$.
If $\xbar_i\subseteq\xbar_m$, then $\alpha_i(\cbar_n\mr{\xbar_i})$ holds because $\cbar_n\mr{\xbar_i}=\abar_n\mr{\xbar_i}$ and $\abar_n$ realizes $\theta(\xbar)$.
On the other hand, if $\xbar_i\cap\xbar_m=\emptyset$, then $\cbar_n\mr{\xbar_i}=\bbar\mr{\xbar_i}=\abar_0\mr{\xbar_i}$ and $\abar_0$ realizes $\theta(\xbar)$.

To finish the proof of the Claim, it suffices to show that for any $j\in U$, $\neg\beta_j(\cbar_n\mr{\xbar_j})$ holds for cofinitely many $n$.
Write $\xbar_j=\xbar'\conc\ybar'$ with $\xbar'=\xbar_j\cap\xbar_m$ and $\ybar'=\xbar_j\setminus\xbar_m$, and
again we split into cases.  If $\ybar'=\emptyset$, i.e., $\xbar_j\subseteq\xbar_m$, then $\cbar_n\mr{\xbar_j}=\abar_n\mr{\xbar_j}$, so $\neg\beta_j(\cbar_n\mr{\xbar_j})$
since $\theta(\abar_n)$.  If $\xbar'=\emptyset$ then $\cbar_n\mr{\xbar_j}=\bbar\mr{\xbar_j}=\abar_0\mr{\xbar_j}$, so again $\neg\beta_j(\cbar_n\mr{\xbar_j})$ since $\theta(\abar_0)$.
Finally, if both $\xbar'$ and $\ybar'$ are non-empty, $\xbar'\conc\ybar'$ is a proper partition of $\xbar_j$.
Let $\bbar'=\bbar\mr{\ybar'}$.  As $\beta_j(\xbar_j)$ is mutually algebraic, choose an integer $s$ such that $\exists^{\le s}\xbar'\beta_j(\xbar',\bbar')$.
As $\{\cbar_n\mr{\xbar'}:n\in\omega\}$ are disjoint, there are at most $s$ $n$'s for which $\beta_j(\cbar_n\mr{\xbar'},\bbar')$ holds, hence
$\neg\beta_j(\cbar_n\mr{\xbar_j})$ holds for cofinitely many $n$.
\endproof

The following Proposition reaps the benefit of a finite, relational language.

\begin{Proposition}  \label{ap}  Suppose $M$ is mutually algebraic in a finite, relational language with every atomic formula having free variables among $\zbar$.
 There is a finite set $\Fcal=\{\phi_i(\xbar_i,\wbar_i):i<m\}$ of quantifier-free $L$-formulas  such that
for every $\xbar\subseteq\zbar$ and every $p\in \S_{\xbar}(M)$ that contains a mutually algebraic formula, there is some $\phi_i\in\Fcal$ and $\ebar_i\in M^{\lg(\wbar_i)}$ such that $\phi_i(\xbar,\ebar_i)\in p$ and is mutually algebraic.
\end{Proposition}

\bp  First, by Theorem~\ref{old}, there is a finite set $\B=\{\delta_i(\xbar_i,\ebar_i):i<n\}$ of mutually algebraic $L(M)$-formulas such that every atomic $R\in L$ is equivalent to a boolean combination of formulas from $\B$.  It follows that for every $\xbar\subseteq\zbar$, every $\gamma(\xbar)\in\F_{\xbar}(M)$ is also equivalent to a boolean combination of $\B$-formulas.
Let $k=\lg(\zbar)$ and let $\B_k$ denote the (finite)  set of $\le k$-conjunctions of formulas from $\B$, i.e., $\B_k=\{\bigwedge\B_0:\B_0\subseteq\B$ and $|\B_0|\le k\}$.
Let
$$\Fcal=\{\phi(\xbar,\wbar):\xbar\subseteq\zbar\ \hbox{and $\phi(\xbar,\hbar)\in\B_k$ for some $\hbar\in M^{\lg(\wbar)}\}$}$$

To see that $\Fcal$ is as desired, fix $\xbar\subseteq\zbar$,  $p\in\S_{\xbar}(M)$ containing a mutually algebraic formula $\gamma(\xbar)$, and a realization $\abar$ of $p$.   
Write $\gamma(\xbar)$ as a disjunction $\bigvee\theta_\ell$, where each $\theta_\ell$ is a conjunction of $\B$-formulas and their negations.
Let $\theta(\xbar)$ be one of the conjuncts for which $\theta(\abar)$ holds.  Then $\theta(\xbar)\in p$ and since $\theta(\xbar)\vdash\gamma(\xbar)$,
$\theta(\xbar)$ is mutually algebraic.  Write $$\theta(\xbar)=\bigwedge_{i\in S}\delta_i(\xbar_i,\ebar_i)\wedge\bigwedge_{j\in U}\neg\delta_j(\xbar_j,\ebar_j)$$
with each $\delta_i(\xbar_i,\ebar_i),\delta_j(\xbar_j,\ebar_j)\in \B$, hence mutually algebraic.  
Apply Lemma~\ref{conj} to $\theta(\xbar)$, obtaining $S_0\in [S]^{\le n}$ as there.  Thus, $\phi(\xbar,\hbar):=\bigwedge_{i\in S_0}\delta_i(\xbar_i,\ebar_i)\in p$ 
and is mutually algebraic.  Visibly, $\phi(\xbar,\hbar)\in\B_k$,
so $\phi(\xbar,\wbar)\in\Fcal$ as required.
\endproof

\end{document}